\documentclass[12pt,reqno]{amsart}

\usepackage{amsmath}
\usepackage{amssymb}
\usepackage{amscd}
\usepackage{pstricks}
\usepackage{color}
\usepackage{mathpazo}
\usepackage{url}
\usepackage{graphicx}
\usepackage{caption}
\usepackage{subcaption}

\usepackage[english]{babel}
\theoremstyle{plain}
\newtheorem{theorem}{Theorem}

\newtheorem{definition}{Definition}

\newtheorem{conjecture}[theorem]{Conjecture}

\newtheorem{remark}[theorem]{Remark}
\newtheorem{question}[theorem]{Question}

\newcommand{\bT}{{\mathbb T}}
\newcommand{\bS}{{\mathbb S}}
\newcommand{\bR}{{\mathbb R}}
\newcommand{\bZ}{{\mathbb Z}}

\newcommand{\mM}{\mathcal{M}}

\newcommand{\mQ}{\mathcal{Q}}

\newcommand{\rmd}{\mathrm{d}}

\begin{document}

\title[Random walks and Lorentz processes]
{Random walks and Lorentz processes}

\author{Domokos Sz\'asz}\footnote{Department of Stochastics,
Budapest University of Technology and Economics, H-1111, Budapest; (szasz at math.bme.hu).}

\maketitle

\vskip3mm
\centerline{\it Dedicated to the memory of P\'al R\'ev\'esz}
\bigskip
\begin{abstract}
Random walks and Lorentz processes serve as fundamental models for Brownian motion. The study of random walks is a favorite object of probability theory whereas that of Lorentz processes belongs to the theory of hyperbolic dynamical systems. Here we first present examples where the method based on the probabilistic approach led to new insights into the study of the Lorentz process. Motivated by a 1981 question of Sinai about limiting laws for planar locally perturbed Lorentz processes, we first derived that - in the plane - local perturbations of homogeneous random walks leave the limit laws and the limiting processes unchanged - independently whether the walk had bounded or unbounded jumps. Afterwards, we obtained probabilistic statements for local perturbations of planar Lorentz processes with finite horizon. (Similar statements for local perturbations of Lorentz processes with infinite horizon is a most interesting open problem!) Often an interesting consequence of a local limit theorem for a process is the recurrence of the same process. In this way our approach also provides an alternative proof for a result of Lenci \cite{MR1992667} about the recurrence of the locally perturbed Lorentz processes with finite horizon. Afterwards an unsolved problem - related to Sinai's 1981 question - is formulated as an analogous problem in the language of random walks.

\end{abstract}
\section{Introduction.}

According to P\'olya's classical theorem \cite{MR1512028}, the simple symmetric random walk (SSRW) on $\mathbb Z^d$ is recurrent if $d=1, 2$ otherwise it is transient. Since the one- and two-dimensional Lorentz processes with a periodic configuration of scatterers share a number of stochastic properties with those of SSRW's, it had also been expected that the analogues of P\'olya's theorem also hold for them. Indeed, in $d = 2$, the recurrence of the periodic finite horizon Lorentz process (FHLP) got settled by ergodic theoretic methods independently by K. Schmidt (\cite{MR1663750}) and J.-P. Conze (\cite{MR1721618}). For the same case, in 2004, with T. Varj\'u, we succeeded to give a probabilistic-dynamical proof in \cite{szasz2004local}. The recurrence of the infinite horizon Lorentz process (IHLP) in the plane was open until 2007, when the method of \cite{szasz2004local} could be extended to the infinite horizon case (see \cite{szasz2007limit}). It is worth mentioning that, in the infinite horizon case, the limit law of the Lorentz process belongs to the non-standard domain of attraction of the normal law in contrast with the finite horizon case where convergence to the normal law (and to the Wiener process) holds with the diffusive scaling.

Sinai asked in 1981 whether probabilistic properties of the  Lorentz process remain valid if one changes the scatterer configuration in a bounded domain. For convergence to the Gaussian law (and to the Wiener process) and, moreover, for the local limit theorem (LLT) we gave an affirmative answer in \cite{MR2527323}.
Of course, other properties are also much important as the speed of correlation decay, first-return and first-hitting times, local times, and first-intersection times, almost sure invariance principle, \dots (cf. \cite{MR2401621}). It is worth mentioning that often the local limit theorem implies the recurrence of the process in question. In this paper, for simplicity, we mostly focus on recurrence properties by remarking that these can be often derived from just ergodic properties of the processes involved.

However, as to Sinai's 1981 question, nothing is known in the infinite horizon case. Since the study of random walks might be much instructive for the behavior of Lorentz processes, we - in Section \ref{sec:RWIJ} - formulate some problems for random walks with unbounded jumps that are most interesting in themselves. Before this closing section, the paper contains a survey of results known for stochastic properties of RW's and of Lorentz processes.

\section{Random walks and Lorentz processes}
\subsection{Random walks}

\begin{definition}[Random walk]
\begin{enumerate}
\item
Let $\{X_n \in \mathbb Z^d| n \ge 0\} $ be independent random variables and for $n \ge 0$ denote
\[
S_n = \sum_{j=1}^n X_j.
\]
Then the Markov chain $S_0, S_1, S_2, \dots, S_n, \dots$ is called a {\rm random walk}.

\item The random walk is called {\rm symmetric} if the $X_n$'s are symmetric random variables.

\item The random walk is {\rm simple} if $|X_n| = 1\, \textrm{holds}\,
 \forall n \ge 0.$

 \item The probabilities
\[
P(X_{n+1}= k| S_n)
\]
are called the {\rm transition (or jump) probabilities} of the random walk.

\item The random walk is called {\rm  translation-invariant} (or classical or homogeneous) if its transition probabilities are  translation-invariant.

\end{enumerate}
\end{definition}

As to basic notions and properties of random walks we refer to the monographs \cite{MR388547,MR3060348,MR2240535}.

\begin{definition}[Locally perturbed random walk]\label{def:LPRW}
Assume $a>0$. If - possibly outside an origo-centered cube $Q_a$ of size $a$ - the translation probabilities of a random walk $S_0, S_1, S_2, \dots, S_n, \dots$ are translation-invariant - , then we say that the random walk is a {\rm locally perturbed random walk} (LPRW) (more precisely an {\rm $a$-locally perturbed random walk}). For simplicity we assume that all transition probabilities are bounded away from $0$.
\end{definition}

\subsection{Sinai billiards and Lorentz processes}
\subsubsection{Sinai billiards}

As far as notations go, we mainly follow \cite{MR2229799} for planar
billiards.

Billiards are defined in Euclidean domains bounded by a finite number
of smooth boundary pieces.  For our purpose a \textit{billiard} is a
dynamical system describing the motion of a point particle in a
connected, compact domain $\mQ \subset \bT^{d} = \bR^{d} /
\bZ^{d}$. In general, the boundary $\partial \mQ$ of the domain is
assumed to be piecewise $C^3$-smooth, i.e.~there are no corner points;
if $0<J < \infty$ is the number of such pieces, we can write $\partial
\mQ = \cup_{1\leq \alpha\leq J}\,\partial \mQ_\alpha$. Connected
components of $ \bT^{d} \setminus \mQ$ are
called \textit{scatterers} and are assumed to be strictly convex. Motion is uniform inside $\mQ$ and specular
reflections take place at the boundary $\partial \mQ$; in other words,
a particle propagates freely until it collides with a scatterer, where
it is reflected elastically, i.e.~following the classical rule that the
angle of incidence be equal to the angle of reflection.

\begin{definition}[Sinai billiard]
A billiard with strictly convex scatterers is called a {\rm Sinai billiard}.
\end{definition}
\begin{remark}
The above notion of a Sinai billiard is more general than its original one where it was supposed that $d=2, J=1$ and $\mQ_1$ was a circle.
\end{remark}

Since the absolute value of the velocity is a first integral of
motion, the phase space of our billiard is defined as the product of
the set of spatial configurations by the $(d-1)$-sphere, $\mM=\mQ\times
\bS_{d-1}$, which is to say that every phase point $x\in \mM$ is of the
form $x=(q,\,v)$, with $q\in \mQ$ and $v\in \bR^{d}$ with norm
$|v|=1$.
According to the reflection rule, $\mM$ is subject to identification
of incoming and outgoing phase points at the boundary
$\partial\mM=\partial\mQ\times \bS_{d-1}$.
The billiard dynamics on $\mM$ is called the \textit{billiard flow}
and denoted by $S^{t}: t \in (-\infty, \infty)$, where $S^{t}: \mM \to
\mM$. The set of points defined by the trajectory going through
$x\in\mM$ is denoted $S^{\bR}x$. The smooth, invariant probability
measure of the billiard flow, $\mu$ on $\mM$, also called the
Liouville measure, is essentially the product of Lebesgue measures on
the respective spaces,
i.e.~$\rmd\mu= {\rm const.}\,  \rmd q \, \rmd v$, where the constant
is $(\textrm{vol} \,\mQ \ \textrm{vol}\, \bS_{d-1})^{-1}$.

The appearance  of collision--free orbits is a distinctive feature of some billiards
  which are said to have infinite horizons.
\begin{definition}[Infinite and finite horizons]\label{def:horizons}

   \begin{enumerate}
  \item
    Denote by $\mM_{\textrm{free}} \subset \mM$ the subset of collision--free
    orbits, i.e.
    \begin{equation*}
      \mM_{\textrm{free}} = \{x \in \mM \,:\, S^{\bR}x
      \cap \partial \mM = \emptyset\}\,.
    \end{equation*}
  \item
    The billiard has \emph{finite horizon} if $\mM_{\textrm{free}} =
    \emptyset$. Otherwise it has \emph{infinite horizon}.
  \end{enumerate}
\end{definition}

\subsubsection{Lorentz processes}
The Lorentz process was introduced in
1905 by H. A. Lorentz \cite{L1905} for the study of a dilute electron
gas in a metal.
While Lorentz considered the motion of a collection of independent
pointlike particles moving uniformly among immovable metalic ions
modeled by elastic spheres, we consider here the uniform motion of a
single pointlike particle in  a fixed array of strictly convex scatterers with
which it interacts via elastic collisions.

Thus defined, the \textit{Lorentz process} is the billiard dynamics of
a point particle on a billiard table $\mQ =  \bR^d
\setminus \cup_{\alpha=1}^\infty \, O_{\alpha}$, where the scatterers
$O_{\alpha}$,  $ 1 \leq \alpha \leq  \infty$, are strictly convex with
$C^3$-smooth boundaries. Generally speaking, it could happen that
$\mQ$ has several connected components. For simplicity,
however, we assume that the scatterers are disjoint and that
$\mQ$ is unbounded and connected.  The phase space of this process is
then given according to the above definition, namely $\mM=\mQ\times
\bS_{d-1}$.

\begin{definition}[Lorentz process]
Select an initial phase point $(q_0, v_0) \in \mM$. Let its trajectory under the Lorentz (or billiard) dynamics be $\{(q_j, v_j)| j \ge 0\}$. Then the process $\{q_j| j \ge 0\}$ is called the {\rm Lorentz process} (cf. Figure 1).
\end{definition}

\begin{definition}[Diffusively scaled flow]
Assume $\{q_n \in \mathbb R^d| n \ge 0\}$ is a random trajectory where, for simplicity,  $d = 2$.
Then its {\rm diffusively scaled flow} $\in C[0,1]\ \  (\textit{or} \ \ \in C[0, \infty])$ is defined as follows:
for $N \in \mathbb Z_+$ denote $W_N(\frac{j}{N}) = \frac{q_j}{\sqrt N}\qquad
(0 \le j \le N\ \  \textit{or} \ \ j \in \mathbb Z_+)$
and define otherwise $L_N(t)( t \in [0,1] \ \ \textit{or} \ \ \mathbb R_+)$ as its piecewise linear, continuous extension.
\end{definition}

We note that, for the Lorentz dynamics,
the Liouville measure $\rmd\mu= \rmd q\, \rmd v$, while invariant, is
infinite. If, however, there exists a regular lattice of rank $d$ for
which we have that, for every point $z$ of this lattice, $\mQ + z =
\mQ$, then we say that the corresponding Lorentz process is
\textit{periodic}. In this case, the Liouville measure is finite (more
exactly, its factor with respect to the lattice is finite). (For simplicity, the lattice of periodicity will be taken as $\mathbb Z^d$.)

\begin{definition}[Periodic Lorentz process]
The Lorentz process is called {\rm periodic} if the configuration $\{O_{\alpha}| 1 \le \alpha < \infty\}$  of its scatterers is $\mathbb Z^d$-periodic.
\end{definition}

\begin{figure}\label{fig:PLP}
\centering
\begin{subfigure}{.45\textwidth}
  \centering
  \includegraphics[width=.87\linewidth]{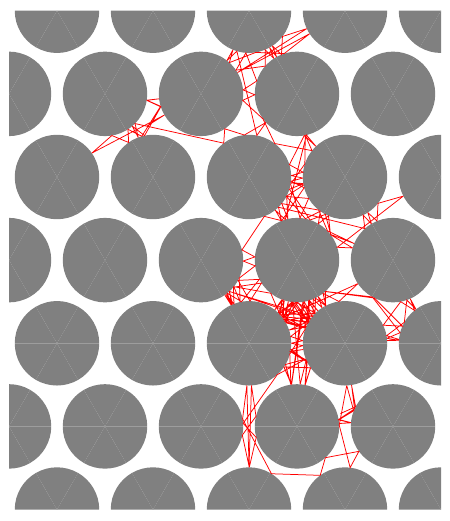}
  \caption{Finite horizon}
  \label{fig:sub1}
\end{subfigure}%
\begin{subfigure}{.5\textwidth}
  \centering
  \includegraphics[width=1.3\linewidth]{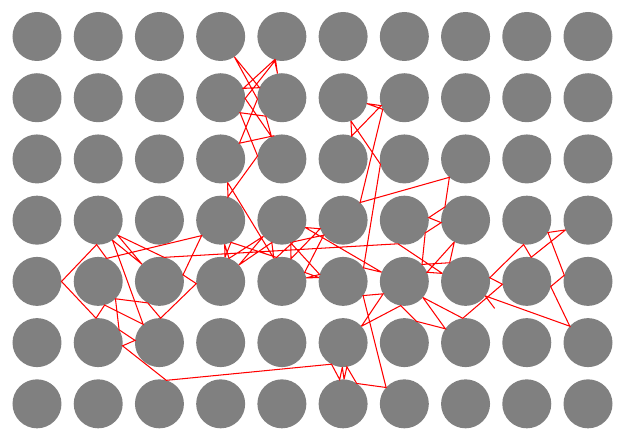}
  \caption{Infinite horizon}
  \label{fig:sub2}
\end{subfigure}
\caption{Periodic Lorentz process}
\label{fig:test}
\end{figure}

Figure 1 shows examples of trajectories (red lines) in periodic Lorentz processes (for the finite vs the infinite horizon cases, resp. cf. Definition \ref{def:horizons}).

\begin{definition}[Locally perturbed (periodic) Lorentz process]
If one changes arbitrarily the scatterer configuration of a periodic Lorentz process in a bounded domain, then we talk about a {\rm locally perturbed Lorentz process (LPLP)}. This process has {\rm finite or infinite horizon}  if the original periodic Lorentz process had finite vs. infinite horizon respectively.
\end{definition}

\subsubsection{Recurrence of stochastic processes}

\begin{definition}[Recurrence]
A stochastic process in $\mathbf Z^d$ or in $\mathbf R^d$ ($d \ge 1$) is {\rm recurrent} if for any bounded subset of $\mathbf Z^d$ (or of $\mathbf R^d$) it is true that the process returns to the subset infinitely often with probability $1$.
\end{definition}
\section{Recurrence of periodic random walks and Lorentz processes in the plane}
In this section we restrict ourselves to the case when - on the one hand - the transition probabilities of the random walk are translation invariant and - on the other hand - the Lorentz process is periodic.
\subsection{Random walks}

Start with the classical theorem of P\'olya.
\begin{theorem}[\cite{MR1512028}]
The SSRW is recurrent if $d=1, 2$ and otherwise it is transient.
\end{theorem}
For more general random walks we refer to results of Chung-Ornstein and Breiman.
\begin{theorem}[\cite{MR133148,MR1163370}]
\begin{enumerate}
\item Assume $(S_0, S_1, S_2, \dots)$ is a translation invariant RW on $\mathbb Z$. If $\mathbf EX_1 = 0$, then the RW is recurrent.
\item Assume $(S_0, S_1, S_2, \dots)$ is a translation invariant RW on $\mathbb Z^2$. If $\mathbf EX_1 = 0$ and $\mathbf E|X_n|^2 < \infty$, then the RW is recurrent.

\end{enumerate}
\end{theorem}

\subsection{Lorentz processes}

Based on the analogy with random walks, for periodic Lorentz processes the exact analogue of P\'olya’s theorem  known for
random walks had been expected.

\subsubsection{Finite horizon}
The first positive result was obtained in \cite{MR789870},
where a slightly weaker form of recurrence was demonstrated: the process almost surely
returns infinitely often to a moderately (actually logarithmically) increasing sequence of
domains. The authors used a probabilistic method combined with the dynamical tools
of Markov approximations. (The weaker form of recurrence was the consequence of the
weaker form of their local limit theorem based on the weaker CLT of \cite{MR606459}.)

An original approach appeared in 1998–99, when independently
Schmidt \cite{MR1663750} and Conze \cite{MR1721618} were, indeed, able to deduce recurrence from the
global central limit theorem (CLT) of \cite{MR1138952} by adding (abstract) ergodic theoretic ideas.

\begin{theorem}[\cite{MR1721618,MR1663750}]
 The planar Lorentz process with a finite horizon is almost surely recurrent.
 \end{theorem}
 Their approach seems, however, to be essentially restricted to the planar finite horizon case.

\subsubsection{Infinite horizon}
For attacking the infinite horizon case, the authors of \cite{szasz2004local} returned to the probabilistic approach via the local limit theorem and first gave a new proof of the theorems of Conze and Schmidt. Finally, in \cite{szasz2007limit}, they could prove a local limit theorem for the Lorentz process in the infinite horizon case that already implied recurrence in this case, too.
\begin{theorem}[\cite{szasz2007limit}]
 The planar Lorentz process with infinite horizon is almost surely recurrent.
\end{theorem}

\section{Stochastic properties of locally perturbed planar Lorentz processes}
\subsection{Finite horizon case}

Answering Sinai's 1981 question we could prove the following
\begin{theorem}[\cite{MR2527323}]
 Consider the diffusively scaled variant  $\{L_N(t)|t \ge 0\}$  of a locally  perturbed finite horizon Lorentz process.
 Then, as $N \to \infty$,  $L_N(t) \Rightarrow W_{\sigma^2}(t)$
(weak convergence in $C[0, \infty]$),
where $W_{\sigma^2}(t)$ is the Brownian Motion with the non-degenerate covariance matrix
$\sigma^2$. The limiting covariance matrix coincides with that for the
unperturbed periodic Lorentz process (cf. \cite{MR1138952}).
 \end{theorem}

We note that the proof of the above theorem uses delicate  recurrence properties of the periodic Lorentz process (cf. \cite{MR2401621}) being interesting in themselves. Actually they are analogues of several properties of classical random walks.
 Moreover, the local version of the above theorem provides an alternative proof for the following recurrence result of Lenci:
 \begin{theorem}[\cite{MR1992667}]
 The locally perturbed planar Lorentz process with a finite horizon is almost surely recurrent.
\end{theorem}

\subsection{Infinite horizon case}
\begin{definition}[Superdiffusively scaled flow]
Assume $\{q_n \in \mathbb R^d| n \ge 0\}$ is a random trajectory where, for simplicity,  $d = 2$.
Then its {\rm superdiffusively scaled flow} $\in C[0,1]\ \  (\textit{or} \ \ \in C[0, \infty])$ is defined as follows:
for $N \in \mathbb Z_+$ denote $L_N(\frac{j}{N}) = \frac{q_j}{\sqrt {N\log n}}\qquad
(0 \le j \le N\ \  \textit{or} \ \ j \in \mathbb Z_+)$
and define otherwise $L_N(t)( t \in [0,1] \ \ \textit{or} \ \ \mathbb R_+)$ as its piecewise linear, continuous extension.
\end{definition}

\begin{conjecture}\label{thm:LPIHLP}
\begin{enumerate}
\item
 Consider the superdiffusively scaled flow  $\{L_N(t)|t \ge 0\}$  for a locally  perturbed infinite horizon Lorentz process.
 Then, as $N \to \infty$,  $L_N(t) \Rightarrow W_{\sigma^2}(t)$
(weak convergence in $C[0, \infty]$),
where $W_{\sigma^2}(t)$ is the Brownian Motion with the non-degenerate covariance matrix
$\sigma^2$. The limiting covariance matrix coincides with that for the
unperturbed periodic Lorentz process (cf. \cite{szasz2007limit}).
\item
 The locally perturbed planar Lorentz process with infinite horizon is almost surely recurrent.

\end{enumerate}
\end{conjecture}

\section{Stochastic properties of symmetric vs. locally perturbed random walks with unbounded jumps}\label{sec:RWIJ}

For understanding the difficulties in proving Conjecture \ref{thm:LPIHLP} it should be instructive to answer the analogous question for LPRW's with unbounded jumps.

For LPRW's with bounded jumps, the first result related to Sinai's question was given in the paper \cite{MR648199}.
\subsection{Reminder on some results of \cite{MR648199}}
We recall a simple example of its main theorem.
\begin{definition}\label{def:LPSRW}
Let $\{S_n| n \ge 0\} $ be a simple RW on $\mathbb Z^2$ such that for $i = 1, 2$
\[
P(X_{n+1}= \pm e_i|S_n) =
\begin{cases}\hskip6mm \frac{1}{4} & \textrm{if}\, \,  S_n\neq (0, 0)\\
{\rm arbitrary}  & {\rm if}\, S_n = (0, 0)
\end{cases}
\]
where $e_1 = (1, 0)$ and $e_2 = (0, 1)$.
\end{definition}
Let, moreover,
\[
U_n(t): = n^{-1/2}S_{[nt]} \hskip1cm t \in [0, 1]
\]
\begin{theorem}[\cite{MR648199}]\label{thm:SZT81}
As $ n \to \infty$
\[
U_n(t) \Rightarrow W(t)
\]
weakly in $C[0, 1]$, where $W$ is the standard planar Wiener process.
\end{theorem}

\begin{remark}\label{thm:LCLTetc}
\begin{enumerate}
\item
By applying the methods of \cite{MR648199,MR2401621}, one can easily see that
\begin{enumerate}
\item
For $U_n(1)$ the global CLT holds;
\item For $U_n(1)$ the local CLT also holds;
\item As a consequence of the later one, the RW of Definition \ref{def:LPSRW} is recurrent.
 \end{enumerate}
\item In fact, Theorem \ref{thm:SZT81} is the special case of a more general theorem of \cite{MR648199} whose statement roughly says that, if one has a CLT for a translation-invariant RW, then changing the jump probabilities in a bounded domain does not change the statement of the CLT. The methods of \cite{MR2401621} also imply the truth of the local limit theorem and recurrence in this generality.
\end{enumerate}
\end{remark}

\subsection{A locally perturbed RW with unbounded jumps}
Following the bounded jumps case, for the unbounded jumps case we will also start with a simple example:
\begin{definition}\label{def:LPRWUJ}
Let  $\{S_n| n \ge 0\} $ be a RW (with unbounded jumps) such that for $i =1, 2$ one has
\[
 \begin{cases}
 P(X_{n+1} = \pm e_i| S_n=(0,0) ) &=1/4\\
 P(X_{n+1} = \pm ne_i| S_n\neq (0,0) )&=
  const. \frac{1}{|n|^3}
 \end{cases}
 \]
 \end{definition}

 Let us explain why we suggest first the study of exactly this example. The simplest example of a periodic Lorentz process is the following:
 all scatterers are circles of radius $R$ with $R_{\rm min} < R < \frac{1}{2}$. This condition ensures that all elements of  $\mM_{\textrm{free}}$ are parallel to one of the axes. In this case the long jumps of the Lorentz process are almost parallel to one of the axes and the distribution of their lengths is asymptotically $ const. \frac{1}{|n|^3}$ and thus belongs to the non-standard domain of attraction of the normal law (cf. \cite{MR1149489,szasz2007limit,D12,MR3231997}).

The  unperturbed RW corresponding to the previous example belongs to the non-standard domain of attraction of the normals law with $\sqrt{n \log n}$ scaling. Denote
\[
V_n(t): = (n \log n)^{-1/2}S_{[nt]} \hskip1cm t \in [0, 1]
\]

A special case of the main result of the work \cite{MR2740407} is the following:
\begin{theorem}[\cite{MR2740407}]\label{thm:PSZ10}
As $ n \to \infty$
\[
V_n(t) \Rightarrow C W(t)
\]
weakly in $C[0, 1]$, where $W$ is the standard planar Wiener process and $C > 0$.
\end{theorem}

\begin{conjecture}\label{Conj}
In the setup of this theorem
\begin{enumerate}
\item
the local version of the limit law for $V_n(1)$ is also true;
 \item the RW defined in Definition \ref{def:LPRWUJ} is recurrent.
 \end{enumerate}
\end{conjecture}

\section{Strongly perturbed RW's}
For locally perturbed random walks, it was sort of expected that a local perturbation should not change the limiting behavior of the RW whether the jumps are bounded or unbounded whatever difficulties the proofs of these statements would bring up. It is, however, an intriguing question under what kind of extended perturbations the classical limiting behavior of the random walk survives. Before a precise formulation of this question let us introduce notations.

\begin{definition}[Strongly perturbed random walks]
Assume $\{0 < a_n|n \ge 1\}$ are such that $ \lim_{n \to \infty} a_n = \infty$. The sequence $S^{a_n}_0, S^{a_n}_1, S^{a_n}_2, \dots, S^{a_n}_n, \dots$ of $a_n$-locally perturbed random walks is called strongly $a_n$-perturbed.
\end{definition}

\begin{question}
For a sequence $\{a_n| \lim_{n \to \infty} a_n = \infty \}$ denote
\[
Z_n(t): = n^{-1/2}S^{a_n}_{[nt]} \hskip1cm t \in [0, 1],
\]
where the jump probabilities of the translation-invariant random walk are as in Definition \ref{def:LPSRW}.
\begin{enumerate}
\item
Find a sequence  $\{a_n| \lim_{n \to \infty} a_n = \infty \}$ and a strongly $a_n$-locally perturbed sequence of random walks (cf. Definition \ref{def:LPRW}) with bounded jumps  such that
\begin{equation}\label{eq:Z}
Z_n(t) \Rightarrow W(t)
\end{equation}
weakly in $C[0, 1]$.

\item
Can you prove Equation (\ref{eq:Z}) for any sequences of $a_n$ with $a_n = o(n^{1/2})$?
\end{enumerate}
\end{question}

\begin{question}
For a sequence $\{a_n| \lim_{n \to \infty} a_n = \infty \}$ denote
\[
Y_n(t): = (n \log n)^{-1/2}S^{a_n}_{[nt]} \hskip1cm t \in [0, 1]
\]
where the jump probabilities of the translation-invariant random walk are as in Definition \ref{def:LPRWUJ}.
\begin{enumerate}
\item
Find a sequence $\{a_n| \lim_{n \to \infty} a_n = \infty \}$  and a strongly $a_n$-locally perturbed sequence of random walks with unbounded jumps  such that
\begin{equation}\label{eq:Y}
Y_n(t) \Rightarrow W(t)
\end{equation}
weakly in $C[0, 1]$.
\item
Can you prove Equation (\ref{eq:Y}) for any sequences of $a_n$ with $a_n = o((n \log n)^{1/2})$?
\end{enumerate}
\end{question}

\section{Acknowledgement} The author expresses his gratitude to Istv\'an Berkes for his precious advice. I am also grateful to the referees for their useful remarks. Support of the NKFIH Research Fund, grant K142169, is thankfully acknowledged. The author thanks the Erwin Schrödinger International Institute for Mathematics and Physics since the questions raised in the paper were initiated during the workshop Rare Events in Dynamical Systems (March 18, 2024 -- March 22, 2024).

\section{Added in proof.} After the closure of this manuscript Marco Lenci has informed me that he can prove statement (2) of Conjecture \ref{thm:LPIHLP}. The author also expresses his gratitude to him for his precious remarks on a previous version of this work.
\bibliographystyle{alpha}
\bibliography{revesz_nodoi_revised}

\end{document}